\documentclass[11pt]{amsart}
\usepackage{amsaddr}
\usepackage{geometry}                
\usepackage{graphicx}
\usepackage{amssymb}
\usepackage{enumerate}
\usepackage{epstopdf}
\usepackage{amssymb}
\usepackage{amsmath}
\usepackage{blkarray} 
\usepackage{comment}
\usepackage{latexsym}
\usepackage{amsthm} 
\usepackage{eucal} 
\usepackage{eufrak}
\usepackage{float}
\usepackage[rightcaption]{sidecap}
\usepackage{graphicx}
\graphicspath{ {images/} }
\usepackage[utf8]{inputenc}
\usepackage[main=english, russian]{babel}
\usepackage{xcolor}
\usepackage{mathrsfs}
\usepackage{array}
\usepackage{framed}
\usepackage[all]{xy}
\usepackage{sidecap}
\usepackage{wrapfig}
\usepackage{multicol}
\usepackage{lscape}
\usepackage{MnSymbol}
\usepackage{mathtools}
\usepackage{cite}
\usepackage{appendix}
\usepackage{comment}
\usepackage{chngcntr}
\usepackage{etoolbox}
\usepackage{lipsum}
\usepackage{rotating}
\usepackage[thinlines]{easytable}
\usepackage{hyperref}
\usepackage{empheq}
\usepackage{lineno}
\DeclareMathAlphabet{\mathpzc}{OT1}{pzc}{m}{it}

\theoremstyle{plain} 
\newtheorem{thm}{Theorem}[section] 
\newtheorem{cor}[thm]{Corollary} 
 
\newtheorem{prop}[thm]{Proposition} 

\theoremstyle{definition} 
\newtheorem{defn}{Definition}

\theoremstyle{remark}

\title[Structural obstruction in chemical reaction networks]{Structural obstruction to the simplicity of the eigenvalue zero in chemical reaction networks}
\author{Nicola Vassena}
\address{Freie Universit{\"a}t Berlin}
\email{nicola.vassena@fu-berlin.de}
\date{\today}

\begin{document}

\begin{abstract}
Multistationarity is the property of a system to exhibit two distinct equilibria (steady-states) under otherwise identical conditions, and it is a phenomenon of recognized importance for biochemical systems. Multistationarity may appear in the parameter space as a consequence of saddle-node bifurcations, which necessarily require a simple eigenvalue zero of the Jacobian, at the bifurcating equilibrium. Matrices with a simple eigenvalue zero are generic in the set of singular matrices: any system whose Jacobian has an algebraically multiple eigenvalue zero can be perturbed to a system whose Jacobian has a simple eigenvalue zero. Thus, one would expect that in applications singular Jacobians are always with a simple eigenvalue zero. However, chemical reaction networks typically consider a fixed network structure, while the freedom rests with the various and different choices of kinetics. Here we present an example of a chemical reaction network, whose Jacobian is either nonsingular or has an algebraically multiple eigenvalue zero. The structural obstruction to the simplicity of the eigenvalue zero is based on the network alone, and it is independent of the value of concentrations and the choice of kinetics. This in particular constitutes an obstruction to standard saddle-node bifurcations.
\end{abstract}

\providecommand{\keywords}[1]
{
  \small	
  \textbf{\textit{Keywords---}} #1
}

\maketitle

\tableofcontents

\section{Introduction}

A \emph{bifurcation} is a sudden qualitative change in the system behavior under a small change in the parameter values. Bifurcation theory is hence a powerful tool to identify parameter areas where certain dynamical behaviors of interest occur. 
We refer to the standard book by Guckenheimer and Holmes \cite{GuHo84} for background on bifurcations. The simplest bifurcations occur according to one parameter, only. Standard and well-studied examples are \emph{saddle-node bifurcations}, that detect \emph{multistationarity}, and \emph{Hopf bifurcations}, that detect \emph{oscillations}. Their spectral condition is a \emph{simple eigenvalue zero} (saddle-node) and \emph{a pair of purely-imaginary complex-conjugated eigenvalues} (Hopf). In a biochemical context, which is the focus of the present paper, saddle-node bifurcations under the assumption of mass action kinetics have been discussed by Conradi et al. \cite{Conradi2007} and Domijan and Kirkilionis \cite{DomKirk09}. Otero-Muras and coauthors used computational methods to detect saddle-node bifurcations in biochemical systems, see for example \cite{Otero18} and the many references therein. For Hopf bifurcation, see the pioneering work by Gatermann et al. \cite{Gat2005}, applying concepts from computer algebra to mass action systems, and Fiedler \cite{F19}, with more general kinetics and in global setting. By proving Hopf bifurcations in more circumstantial systems, Conradi et al. detected oscillations in a mixed-mechanism phosphorylation system \cite{conradietal19}, Boros and Hofbauer in planar deficiency-one mass-action systems \cite{Boros21}, Hell and Rendall in the MAP kinase cascade \cite{Hell16}.\\

In the present paper, we investigate the necessary spectral condition for saddle-node bifurcations: a simple eigenvalue zero of the Jacobian, at the equilibrium. Saddle-node bifurcations are often invoked in the quest of finding bistable systems: a typical bistable scheme considers two connecting saddle-node bifurcations, that give rise to bistability and hysteresis phenomena \cite{AnFeSo04}. More in general, a saddle-node bifurcation occurs, e.g., when two equilibria, one stable and one unstable, collide in a single saddle equilibrium and disappear. Hence, such bifurcations point at parameter regions where multistationarity occurs. Multistationarity is the property of a chemical system to exhibit two or more distinct equilibria, co-existing under otherwise identical conditions, and the phenomenon has been proposed as an explanation for many epigenetic processes, including cell differentiation \cite{ThomKauf01}. As a consequence of its importance, multistationarity in chemical reaction networks has been extensively studied, via different methods. See among others the works by Rendall and coauthors \cite{rendall1, rendall2, rendall3}, Dickenstein et al. \cite{Dickensteinetal19}, Shiu and de Wolff \cite{Shiu19}, Feliu et al. \cite{Timo20},  in a mass action context; Soul{\'e}\cite{Soule2003}, Craciun and Feinberg \cite{Crafei05, Crafei06}, Mincheva and Roussel \cite{MR07}, Banaji and Craciun \cite{BaCra10}, Joshi and Shiu \cite{JoShiu13}, Banaji and Pantea \cite{BaPa16}, Conradi et al. \cite{CFMW17}, in more generality.\\

Finding equilibrium bifurcations for a parametric vector field $f(x,\lambda)$ means solving a system of constraints. First, $f(x,\lambda)$ must satisfy the \emph{equilibrium constraints} $f(\bar{x}, \bar{\lambda})=0$, at a certain $\bar{x}$ and $\bar{\lambda}$. Second, \emph{bifurcation conditions} on the Jacobian $\partial_x f(\bar{x},\bar{\lambda})$ must be satisfied. Bifurcation conditions typically comprise the necessary spectral condition on the eigenvalues of the Jacobian and further sufficient nondegeneracy conditions involving higher-order derivatives. Solving such a system of constraints may be very demanding in parametric systems. The point of view of \emph{genericity theory} advocates looking only for the necessary spectral condition, as this is sufficient in ``most" applications to conclude a bifurcation result. This approach brings an obvious advantage in simplifying the mathematical analysis. In fact, bifurcation theory has been historically developed in a general genericity framework \cite{soto73}. A property of a set is generic if it holds on an open and dense subset. In particular, a saddle-node bifurcation happens generically in one-parameter families of vector fields, at an equilibrium with a singular Jacobian. Thus, in the words of Guckenheimer and Holmes \cite[p.~149]{GuHo84}, \textit{ one expects that the zero eigenvalue bifurcations encountered in applications will be saddle-nodes.} This expectation comes from the intrinsic parameter uncertainties in experiments: only generic features are detectable. Interestingly and quite ambiguously, \cite{GuHo84} continues as follows. \textit{If they are not, then there is probably something special about the formulation of the problem which restricts the context so as to prevent the saddle-node from occurring.}\\

Chemical reaction networks theory considers systems of \emph{Ordinary Differential Equations} (ODEs) that are built on two elements: a network structure and parametric reaction rates (kinetics). The network structure is typically considered as given, and it fixes which reactions take place, that is, which reactants react to which products. The parametric nonlinearities model the mathematical laws of the reactions, and they can be chosen with quite a freedom. Standard parametric families of nonlinearities are often considered as reaction rates. For instance, mass action \cite{HFJ72} kinetics, Michaelis-Menten \cite{MM13} kinetics, Hill's kinetics \cite{Hill10}, are important classes of nonlinearities in this context. These relevant kinetics can be grouped and generalized in the definition \ref{chemfeas} below, which assumes that the reaction rate $r_j$ of a reaction $j$ is a positive monotone function of the concentrations of its reactants. The motivating observation for the present paper is that the network structure is fixed, while the nonlinearities are relatively free. Hence, the genericity viewpoint may apply due to the freedom of the reaction rates, or may not apply due to the fixed structure of the network. Specifically, this paper investigates whether the network structure, alone, might be ``\emph{something special about the formulation of the problem}'' preventing a saddle-node bifurcation to occur.\\

We present an example that answers two questions, both affirmatively. The first, quite general and qualitative, is
$$\text{\textbf{Q1:}$\quad$Can the network structure alone prevent generic properties?}$$
Here, ``generic'' indicates properties that are generic in more general families of vector fields, and consequently might be expected to be generic also in the network setting. In detail, we answer the following question.
$$\text{\textbf{Q2:}$\quad$Can the network structure alone prevent the simplicity of the eigenvalue zero?}$$
We present an example of a network for which the associated ODEs system admits a singular Jacobian, but never a Jacobian with a simple eigenvalue zero, for any choice of the reaction rates. This in particular forbids the application of standard saddle-node bifurcation theorems. In general, such an example shows that genericity methods must be handled with extreme care when dealing with chemical reaction networks.\\

The paper is structured as follows. Section \ref{crn} introduces the setting of chemical reaction networks theory. In section \ref{almult} we preliminary provide linear algebra conditions for the simplicity of the eigenvalue zero of general square matrices. In section \ref{mainex} we present our main example of a chemical reaction network, whose associated Jacobian is either nonsingular or possesses an algebraically multiple eigenvalue zero. Sections \ref{almult} and \ref{mainex} are in principle self-contained. For a more in-detail explanation of the intuition behind our example, however, section \ref{CSPCS} introduces the network tools used to design such an example: Child Selections and Partial Child Selections. Section \ref{mainexrev} revisits our main example, explaining it in such network terms. Section \ref{discussion} concludes with the discussion.\\

\textbf{Acknowledgments:} This work has been supported by the Collaborative Research Center (SFB) 910 of the Deutsche Forschungsgemeinschaft (DFG, German Research Foundation)-project number 163436311-SFB 910. We thank Jia-Yuan Dai for his encouragement in working on this problem.

\section{Chemical Reaction Networks}\label{crn}

We briefly present here the setting of Chemical Reaction Networks. A chemical reaction network $\mathbf{\Gamma}$ is a pair of sets $\{\mathbf{M},\mathbf{E}\}$: $\mathbf{M}$ is the set of species or chemicals or metabolites, and $\mathbf{E}$ is the set of reactions. Both sets are finite with cardinalities $|\mathbf{M}|=M$ and $|\mathbf{E}|=E$. Letters $m, n \in \mathbf{M}$ and $j \in \mathbf{E}$ refer to species and reactions, respectively.\\

A reaction $j$ is an ordered association between two positive linear combinations of species:
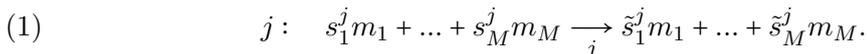
\begin{equation} \label{reactionj}
 j: \quad s^{j}_1m_1+...+s^{j}_Mm_M \underset{j}{\longrightarrow} \tilde{s}^{j}_1m_1+...+\tilde{s}^{j}_Mm_M.
\end{equation}
The nonnegative coefficients $s^{j},\tilde{s}^j$ are called \emph{stoichiometric} coefficients. Chemical networks deal with integer stoichiometric coefficients and typically 0, 1, or 2. However, we can freely consider real $s^{j}_m,\tilde{s}^j_m \in {\mathbb{R}_{\ge 0}}$, as we have no mathematical reason for any restriction. Species appearing at the left (right) hand side of \eqref{reactionj} with nonzero coefficient are called \emph{reactants} (\emph{products}) of reaction $j$. Many chemical systems are open systems: they exchange chemicals with the outside environment. \emph{Inflow reactions} are then reactions with no reactants ($s^j=0$) and \emph{outflow reactions} are reactions with no products ($\tilde{s}^j=0$).\\

The $M \times E$ stoichiometric matrix $S$ is the matrix of all ordered stoichiometric coefficients:
\begin{equation}\label{smatrix}
S_{mj}:= \tilde{s}^j_{m}-s^j_m,
\end{equation}
where $\tilde{s}^j_{m}$ is the stoichiometric coefficient of $m$ as product of $j$, and $s^j_m$ is the stoichiometric coefficient of $m$ as reactant of $j$. With this construction, a fixed order is assigned to each reaction. In particular, we model a reversible reaction
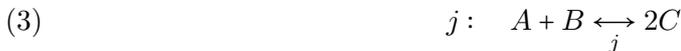
\begin{equation}
j: \quad A+B \underset{j}{\longleftrightarrow} 2C
\end{equation}
simply as two irreversible reactions
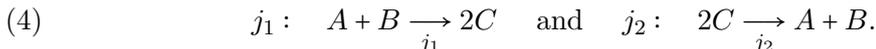
\begin{equation} \label{j_1}
j_1: \quad A+B \underset{j_1}{\longrightarrow} 2C\quad \text{ and }\quad j_2: \quad 2C \underset{j_2}{\longrightarrow} A+B.
\end{equation}

We use the notation $S^j$ to refer to the column of the stoichiometric matrix $S$ associated to the reaction $j$. For example, in a network of four species $\{A,B,C,D\}$, reaction $j_1$ in \eqref{j_1} is represented as the $j_1^{th}$ column of the stoichiometric matrix $S$ as
\begin{equation}
S^{j_1}=
\begin{blockarray}{cc}
 & j_1 \\
\begin{block}{c(c)}
  A & -1\\
  B &  -1\\
  C & 2\\
  D & 0\\
\end{block}
\end{blockarray}\;.
\end{equation}
Note that stoichiometric columns associated with outflow (inflow) reactions always have only negative (positive) entries.\\

Let $x \ge 0$ be the $M$-vector of species concentrations. Under the assumption that the reactor is well-mixed, spatially homogeneous and isothermal, the dynamics $x(t)$ of the concentrations satisfy the following system of ordinary differential equations:
\begin{equation} \label{ODE}
\dot{x}=g(x):=S \mathbf{r}(x),
\end{equation}
where $S$ is the $M \times E$ stoichiometric matrix \eqref{smatrix} and $\mathbf{r}(x)$ is the $E$-vector of the \emph{reaction functions} (kinetics). Without any reactant, we consider as constant the reaction function of inflow reactions $j_f$:
\begin{equation}
r_{j_f}(x)\equiv F_{j_f}.
\end{equation}
For any other reaction $j$, we do not impose any specific form of such functions, requiring only that $r_j$ is \emph{monotone chemical}, as defined in the following definition.
\begin{defn}[monotone chemical function] \label{chemfeas} 
A function $r_j$ is \emph{monotone chemical} if 
\begin{enumerate}
\item $r_j$ depends only on the concentrations of the reactants to reaction $j$:
$$\frac{\partial r_j (x)}{\partial x_m} 
\neq 0 \text{ if, and only if, $m$ is a reactant of $j$}.$$
\item $r_j$ is nonnegative:
$$r(x)\ge0,\text{ for every $x\ge0$},$$
with $r_j(x)=0$ if, and only if, $x_m=0$, for some $m$ reactant of $j$.
\item $\frac{\partial r_j (x)}{\partial x_m} > 0$, for any $m$ reactant of reaction $j$.
\end{enumerate}
\end{defn}

Definition \ref{chemfeas} is standard in many mathematical contributions on chemical networks. Mass action \cite{HFJ72}, Michaelis-Menten \cite{MM13}, Hill's kinetics \cite{Hill10} are important reaction schemes with a wide range of mathematical and biological applications, and they all follow definition \ref{chemfeas}. We use the notation 
$$0 < r_{jm}:=\frac{\partial}{\partial x_m} r_j(x),$$
to refer to the strictly positive partial derivatives. Then, the $E \times M$ \emph{reactivity matrix} $R$ is defined as
\begin{equation} \label{Req}
R_{jm}:=\frac{\partial}{\partial x_m} r_j(x) =
\begin{cases}
r_{jm}\quad \text{if $m$ is a reactant of $j$}\\
0\quad\quad\text{otherwise}
\end{cases}.
\end{equation}
Finally, we define the Jacobian $G$ of the system \eqref{ODE},
$$G := g_x (x) = S R,$$
where $S$ is the stoichiometric matrix and $R$ is the reactivity matrix defined above.\\

\section{Algebraic multiplicity of the eigenvalue zero} \label{almult}
We study the \emph{adjugate matrix} of the Jacobian matrix $G$ to address the geometric and algebraic multiplicity of the eigenvalue zero of $G$. We gather some linear algebra facts \cite{HJ13}. For any $M \times M$ matrix $B$, the adjugate matrix $\operatorname{Adj}({B})$ is the transpose of its cofactor matrix. $\operatorname{Adj}({B})$ satisfies
\begin{equation}
{B} \operatorname{Adj}({B})=\operatorname{Adj}({B}) {B}=\operatorname{det}{B} \; \operatorname{Id}_M,
\end{equation}
where $\operatorname{Id}_M$ is the $M$-dimensional identity matrix. In particular, we have the following straightforward characterizations:
\begin{enumerate}
\item $\operatorname{Adj}(B)$ is invertible if, and only if, $B$ is invertible;
\item $\operatorname{Rank}(B)=M-1$ if, and only if $\operatorname{Rank}(\operatorname{Adj}(B))=1$;
\item $\operatorname{Rank}(B) \le M-2$  if, and only if $\operatorname{Adj}(B)=0$.
\end{enumerate}
Thus $B$ has a geometrically simple eigenvalue 0 if, and only if,
\begin{equation}
\begin{cases}
\operatorname{det} B = 0;\\
\operatorname{Adj}(B) \neq 0.\\
\end{cases}
\end{equation}
Let $\mu_1, ..., \mu_M$ be the eigenvalues of $B$ counted with the respective algebraic multiplicity. Assume $B$ is nonsingular. Then
\begin{equation}\label{tradeig}
\begin{split}
\operatorname{tr}\operatorname{Adj} (B) &=  \operatorname{det}B \operatorname{tr} B^{-1}\\
&=    \prod_{m} \mu_m \; \sum_{m}\frac{1}{\mu_m}\\ 
&=\sum_{i=1}^M \prod_{m\neq i} \mu_m.
\end{split}
\end{equation}
Clearly, \eqref{tradeig} extends to singular matrices $B$, by continuity. In this latter case, we conclude that the algebraic multiplicity of the eigenvalue 0 is exactly one if, and only if, 
$$\operatorname{tr}\operatorname{Adj}(B)\neq 0.$$
In conclusion, we have proved the following proposition.
\begin{prop}\label{mainprop}
The Jacobian $G$ possesses an algebraically simple eigenvalue zero if, and only if, 
\begin{equation}
\begin{cases}
\operatorname{det}G=0;\\
\operatorname{tr}\operatorname{Adj}(G) \neq 0.
\end{cases}
\end{equation}
\end{prop}

As a corollary of interest, we have the following.
\begin{cor}\label{maincor}
Consider the system \eqref{ODE}, and its Jacobian $G$. Assume that, for any choice of $(\mathbf{r}, x)$ such that $\operatorname{det}G=0$, we have that $\operatorname{tr}\operatorname{Adj}(G) = 0$.
Then either
\begin{enumerate}[(i)]
\item $G$ is nonsingular
\end{enumerate}
or
\begin{enumerate}[(i)]
 \setcounter{enumi}{1}
\item $G$ has an eigenvalue zero of algebraic multiplicity strictly bigger than 1.
\end{enumerate}
\end{cor}

\proof[Proof of Corollary \ref{maincor}]
If $\operatorname{det}G \neq 0$, then $G$ is nonsingular. Else, if $\operatorname{det}G=0$, we have by assumption that $\operatorname{tr}\operatorname{Adj}(G) = 0$. By Proposition \ref{mainprop}, we have that the eigenvalue zero is not algebraically simple.
\endproof

The next section presents an example where the assumptions of Corollary \ref{maincor} hold.

\section{Main Example} \label{mainex}

We present a network whose associated system of ODEs admits a singular Jacobian, but never a Jacobian with a simple eigenvalue zero, for any choice of monotone chemical functions endowing the network. In particular, the algebraic multiplicity of the eigenvalue zero is always either $0$ or $>1$. The network possesses 4 species and 6 reactions, and it is as follows.
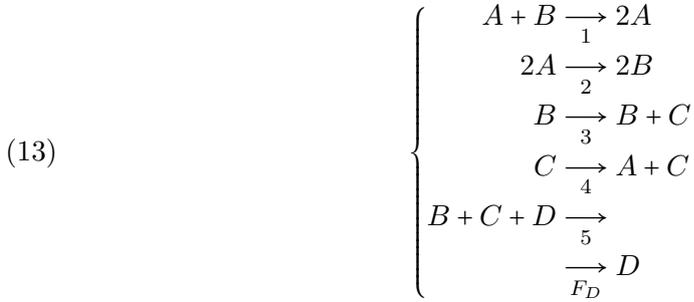
\begin{equation}
\begin{cases}
\begin{split}
A+B &\underset{1}{\longrightarrow} 2A\\
2A &\underset{2}{\longrightarrow} 2B\\
B &\underset{3}{\longrightarrow} B+C\\
C &\underset{4}{\longrightarrow} A+C\\
B+C+D &\underset{5}{\longrightarrow}\\
&\underset{F_D}{\longrightarrow} D\\
\end{split}
\end{cases}
\end{equation}
Reaction $5$ is an outflow from $B, C, D$, while $F_D$ is an inflow to $D$. The system of ODEs is then
\begin{equation}
\dot{x}= S \mathbf{r}(x) =
\begin{blockarray}{ccccccc}
  & 1 & 2 & 3 & 4 & 5 & F_D\\
\begin{block}{c(cccccc)}
 A & 1 & -2 & 0 & 1 & 0 & 0 \\
 B & -1 & 2 & 0 & 0 & -1 & 0 \\
 C & 0  & 0 & 1 & 0 & -1 & 0\\
 D & 0 & 0 & 0 & 0 & -1 & 1 \\
\end{block}
\end{blockarray}
\begin{pmatrix}
r_1(x_A,x_B)\\r_2(x_A)\\r_3(x_B)\\r_4(x_C)\\r_5(x_B,x_C,x_D)\\F_D
\end{pmatrix}
 \end{equation}
We have added labels to rows and columns of the stoichiometric matrix $S$ for simplicity of reading.
In expanded form, the system reads:
\begin{equation}\label{mainexsyst}
\begin{cases}
\dot{x}_A=r_1(x_A,x_B)-2r_2(x_A) + r_4(x_C)\\
\dot{x}_B=-r_1(x_A,x_B)+2r_2(x_A)-r_5(x_B,x_C,x_D)\\
\dot{x}_C= r_3(x_B)-r_5(x_B,x_C,x_D)\\
\dot{x}_D=-r_5(x_B,x_C,x_D)+F_D\\
\end{cases}
\end{equation}

We carry out the analysis in a purely algebraic manner: we investigate the structural relation between $\operatorname{det}G$ and $\operatorname{tr}\operatorname{Adj}(G)$, independently of the chosen value of $x$. In particular, we do not even a priori require that $x$ is an equilibrium. However, preliminarly note that the network does admit an equilibrium: the vector
$$\mathbf{r}_0:=(r_1,r_2,r_3,r_4,r_5,F_D)^T=(r,r,r,r,r,r)^T,$$
$r \in \mathbb{R}_{>0}$, is a positive right kernel vector of the stoichiometric matrix, i.e., it solves
$$0=S\mathbf{r}_0.$$

The Jacobian matrix $G$ of the system reads:
\begin{equation}G=
\begin{pmatrix}
  r_{1A} -2r_{2A} & r_{1B} & r_{4C} & 0\\
    - r_{1A} +2r_{2A} & -r_{1B}-r_{5B} & -r_{5C}& -r_{5D}\\
   0 & r_{3B}-r_{5B} & -r_{5C} &  -r_{5D}\\
0&-r_{5B} & -r_{5C}& -r_{5D}\\
\end{pmatrix}
 \end{equation}
We verify whether the matrix $G$ admits a simple eigenvalue zero. According to section \ref{almult}, we check whether there is a positive solution to the system:
\begin{equation}\label{amsystex}
\begin{cases}
\operatorname{det}G=0\\
\operatorname{tr}\operatorname{Adj}(G)\neq 0\\
\end{cases},
\end{equation}
where $\operatorname{Adj}(G)$ is the adjugate matrix of $G$. We compute the deteminant of $G$ obtaining
\begin{equation}
\operatorname{det}G=(r_{1A}-2r_{2A})r_{3B}r_{4C}r_{5D}.
\end{equation}
Hence, $G$ is singular if, and only if, $\lambda:=r_{1A}-2r_{2A}=0$. We compute the adjugate of $G$:
\begin{small}
\begin{equation}\operatorname{Adj}(G)=
\begin{pmatrix}
 0 & -r_{3B}r_{4C}r_{5D} & -r_{1B}r_{4C}r_{5D} & (r_{1B}+r_{3B})r_{4C}r_{5D}\\
  0 & 0 & \lambda \; r_{4C}r_{5D} &  -\lambda \; r_{4C}r_{5D} \\
    \lambda \; r_{3B}r_{5D}  &  \lambda \; r_{3B}r_{5D}  & 0 &   -\lambda \; r_{3B}r_{5D} \\
- \lambda \; r_{3B}r_{5C} & -\lambda \; r_{3B}r_{5C}  &  -\lambda \; r_{5B}r_{4C} &  \lambda (r_{5B}r_{4C} +r_{3B}r_{5C} -r_{3B}r_{4C}) \\
\end{pmatrix}.
 \end{equation}
 \end{small}
 
 Clearly, the trace of $\operatorname{Adj}(G)$ reads
 \begin{equation}
  \operatorname{tr}\operatorname{Adj(G)}=\lambda (r_{5B}r_{4C} +r_{3B}r_{5C} -r_{3B}r_{4C}),
 \end{equation}
and $\operatorname{tr}\operatorname{Adj(G)}=0$ whenever $\lambda=0$. Thus, system \eqref{amsystex} is never satisfied: the algebraic multiplicity of the eigenvalue zero is strictly bigger than 1, for any choice of $\mathbf{r}$ and $x$ for which $G$ is singular. In particular, note that for a singular $G$, the adjugate matrix reads:

\begin{equation}\operatorname{Adj}\big{|}_{\operatorname{det}G=0}(G)=
\begin{pmatrix}
 0 & -r_{3B}r_{4C}r_{5D} & -r_{1B}r_{4C}r_{5D} & (r_{1B}+r_{3B})r_{4C}r_{5D}\\
  0 & 0 & 0 &  0 \\
  0  &  0& 0 & 0 \\
0 &0 & 0 &  0\\
\end{pmatrix}.
 \end{equation}
Hence, $\operatorname{Adj}(G)\big{|}_{\operatorname{det}G=0} \neq 0$, and the Jacobian $G$ has a geometrically simple but algebraically multiple eigenvalue zero. The precise algebraic multiplicity cannot be asserted from the adjugate matrix alone. An explicit symbolic computation of the eigenvalues shows that the Jacobian $G$, when singular, always possesses an algebraically double eigenvalue zero.

\section{Child Selections and Partial Child Selections} \label{CSPCS}

We introduce the main tools enabling us to discuss the problem in network language. This provides clarification on the design of our main example of section \ref{mainex}.
\begin{defn}[Child Selections] \label{CSDEF}
A \emph{Child Selection} is an injective map 
$$\mathbf{J}: \textbf{M} \longrightarrow \textbf{E},$$ 
which associates to every species $m \in \textbf{M}$ a reaction $j \in \textbf{E}$ such that $m$ is a reactant of reaction $j$.
\end{defn}
Let $S^\mathbf{J}$ indicate the matrix whose $m^{th}$ column is the $\mathbf{J}(m)^{th}$ column of $S$. In particular, the columns of $S^\mathbf{J}$ correspond one-to-one, and following the order, to the reactions
$$\mathbf{J}(m_1), \; \mathbf{J}(m_2), \;...\; , \;\mathbf{J}(m_{M-1}), \;\mathbf{J}(m_{M}).$$
We associate to each Child Selection $\mathbf{J}$ the coefficient
$$\alpha_\mathbf{J} : = \operatorname{det} S^\mathbf{J}.$$
The Jacobian determinant of $G$ can be expressed in terms of Child Selections:
\begin{prop} \label{vita}
Let $G$ be a network Jacobian matrix, in the above settings. Then:
\begin{equation}
\operatorname{det}G=\sum_\mathbf{J} \alpha_\mathbf{J} \; \prod_{m \in \mathbf{M}} r_{\mathbf{J}(m)m},
\end{equation}
\end{prop}
The sum runs on all Child Selections. For a proof, see  \cite{VGB20}. We call the coefficient $\alpha_\mathbf{J}$ \emph{behavior coefficient}. Depending on the sign of $\alpha_{\mathbf{J}}$ we classify a Child Selection as follows. We call a Child Selection $\mathbf{J}$ \emph{zero} if $\alpha_\mathbf{J}= 0$. On the contrary, we call $\mathbf{J}$ a \emph{nonzero} Child Selection if $\alpha_\mathbf{J}\neq0$.\\

We turn now to a related concept: the \textit{Partial Child Selections} (PCS), complimentarily useful to analyze  $\operatorname{Adj}(G)$.

\begin{defn} [Partial Child Selections]
A \emph{Partial Child Selection} $\mathbf{J}^{\vee m_i}$ is an injective map
$$\mathbf{J}^{\vee m_i} : \mathbf{M}\setminus \{m_i\} \longrightarrow \mathbf{E},$$
which associates to each species $m \neq m_i$ a reaction $j$ such that $m$ is a reactant of $j$.\\
\end{defn}

Without loss of generality, assume $1, ... , i, ... , M$. In analogy to the submatrix $S^{\mathbf{J}}$ for a Child Selection $\mathbf{J}$, the expression $S^{\mathbf{J}^{\vee m_i}}$ indicates the $M \times (M-1)$ matrix with columns corresponding one-to-one, and following the order, to the reactions
$$\mathbf{J}^{\vee m_i}(m_1),\;...\;,\;\mathbf{J}^{\vee m_i}(m_{i-1}), \;\mathbf{J}^{\vee m_i}(m_{i+1}),\;...\;,\; \mathbf{J}^{\vee m_i}(m_M).$$
That is, the first column is the stoichiometric column $S^{j_1}$ of the reaction $j_1=\mathbf{J}^{\vee m_i}(m_1)$. Analogously, the $i^{th}$ column is the stoichiometric column $S^{j_i}$ of the reaction $j_i=\mathbf{J}^{\vee m_i}(m_{i+1})$, and so on. We associate to each Partial Child Selection $\mathbf{J}^{\vee m_i}$ the behavior coefficient:
$$\beta_{\mathbf{J}^{\vee m_i}} : =  \operatorname{det}  S^{\mathbf{J}^{\vee m_i}}_{\vee m_i},$$
where the notation $S^{\mathbf{J}^{\vee m_i}}_{\vee m_i}$ indicates the $(M-1)\times(M-1)$ matrix obtained from $S^{\mathbf{J}^{\vee m_i}}$ by removing the $m_i^{th}$ row. If the behavior coefficient $\beta_{\mathbf{J}^{\vee m_i}}$ is zero (nonzero) we call the Partial Child Selection $\mathbf{J}^{\vee m_i}$ \emph{zero} (\emph{nonzero}), accordingly.\\

The role of Partial Child Selections in the analysis of the present paper is clarified by the following proposition.
\begin{prop}\label{AdG}
Let $G$ be the Jacobian matrix of the system \eqref{ODE} and let $\operatorname{Adj}(G)^m_m$ indicate the $m^{th}$ diagonal entry of its adjugate. Then the following expansion holds:
\begin{equation}
\operatorname{Adj}(G)^m_m= \sum_{\mathbf{J}^{\vee m}} \; \beta_{\mathbf{J}^{\vee m}} \;\prod_{n \neq m} r_{\mathbf{J}^{\vee m}(n)n},
\end{equation}
where $\mathbf{J}^{\vee m}$ are Partial Child Selections. The sum runs on all Partial Child Selections $\mathbf{J}^{\vee m}$. In particular, the trace of the adjugate can be expanded as
\begin{equation}
\operatorname{tr}\operatorname{Adj}(G)=\sum_{m\in \mathbf{M}} \sum_{\mathbf{J}^{\vee m}} \; \beta_{\mathbf{J}^{\vee m}} \;\prod_{n \neq m} r_{\mathbf{J}^{\vee m}(n)n}.
\end{equation}
\end{prop}

\proof[Proof of Proposition \ref{AdG}]
The $m^{th}$ diagonal entry of the adjugate matrix of $G$ is
\begin{equation}
\operatorname{Adj}(G)^m_m = \operatorname{det}(G^{\vee {m}}_{{\vee m}}),
\end{equation}
where $G^{\vee {m}}_{{\vee m}}$ indicates the cofactor of $G$ obtained removing the $m^{th}$ row and the $m^{th}$ column. 
We recall that $G=SR$, where $S$ is the stoichiometric matrix (def.\eqref{smatrix}) and $R$ is the reactivity matrix (def.\eqref{Req}). We analyze the expression $\operatorname{Adj}(G)^m_m$ using Cauchy-Binet formula.
\begin{equation}
\begin{split}
\operatorname{Adj}(G)^m_m=& \operatorname{det}(G^{\vee {m}}_{{\vee m}})\\
=&\operatorname{det}(S_{\vee {m}}R^{{\vee m}})\\
=&\sum \limits_{\mathcal{E} \in \mathcal{E}^{M-1}} \operatorname{det} S^{\mathcal{E}}_{\vee m} \; \operatorname{det} R_{\mathcal{E}}^{\vee m}.
\end{split} 
\end{equation}

Note that $\operatorname{det}R_{\mathcal{E}}^{\vee m}\neq 0$ if, and only if, there exists a Partial Child Selection $\mathbf{J}^{\vee m}$, such that $\mathbf{J}^{\vee m}(\mathbf{M} \setminus \{m\})= \mathcal{E}$. Then,

\begin{equation}
\begin{split}
\sum \limits_{\mathcal{E} \in \mathcal{E}^{M-1}} \operatorname{det} S^{\mathcal{E}}_{\vee m} \; \operatorname{det} R_{\mathcal{E}}^{\vee m}&= \sum \limits_{\mathcal{E}=\mathbf{J^{\vee m}}(\mathbf{M} \smallsetminus \{m\})} \operatorname{det}S_{\vee m}^{\mathcal{E}} \; \operatorname{det} R_{\mathcal{E}}^{\vee m} \\
&=\sum \limits_{\mathcal{E}=\mathbf{J^{\vee m}}(\mathbf{M} \smallsetminus \{m\})} \operatorname{det}S_{\vee m}^{\mathcal{E}} \; \operatorname{sgn}(\mathbf{J}^{\vee m}) \; \prod_{n \neq m} r_{\mathbf{J}^{\vee m}(n)n}
\end{split}
\end{equation}
Above, we have expanded $\operatorname{det} R_{\mathcal{E}}^{\vee m}$ via Leibniz formula. The expression $\operatorname{sgn}(\mathbf{J}^{\vee m})$ indicates the signature (or parity) of $\mathbf{J}^{\vee m}$. Then,  \begin{equation}
\operatorname{sgn}(\mathbf{J}^{\vee m}) \, \operatorname{det}S^{ \mathcal{E} =\mathbf{J^{\vee m}}(\mathbf{M} \smallsetminus \{m\})}=  \operatorname{det} S^{{\mathbf{J^{\vee m}}}}_{\vee m}=\beta_{\mathbf{J^{\vee m}}}, 
\end{equation}
which leads to the desired equality: 
\begin{equation} \label{nextproof}
\operatorname{Adj}(G)^m_m= \sum \limits_{\mathbf{J}^{\vee m}}\beta_{\mathbf{J^{\vee m}}} \; \prod_{n \neq m} r_{\mathbf{J}^{\vee m}(n)n}.
\end{equation}
Moreover,
\begin{equation}
\operatorname{tr}\operatorname{Adj}(G)=\sum_m \operatorname{Adj}(G)^m_m = \sum_m  \sum \limits_{\mathbf{J^{\vee m}}}\beta_{\mathbf{J^{\vee m}}} \; \prod_{n \neq m} r_{\mathbf{J}^{\vee m}(n)n}.
\end{equation}
\endproof

\section{Main example revisited}\label{mainexrev}

In this section we explain the design of our main example of section \ref{mainex} in light of the tools introduced in section \ref{CSPCS}. The system \eqref{mainexsyst} possesses only three Child Selections, 
$$\mathbf{J}_1:=\mathbf{J}_1(A,B,C,D)=(1,3,4,5),$$
$$\mathbf{J}_2:=\mathbf{J}_2(A,B,C,D)=(2,3,4,5),$$
and
$$\mathbf{J}_3:=\mathbf{J}_1(A,B,C,D)=(2,1,4,5).$$
The first observation is that the two stoichiometric columns $S^{\mathbf{J}_1(A)}$ and $S^{\mathbf{J}_2(A)}$ are linearly dependent. In particular, $S^{\mathbf{J}_1(A)}=-2\;S^{\mathbf{J}_2(A)}$. This implies that $\mathbf{J}_3$ is a zero Child Selection since
$$\alpha_{\mathbf{J}_3}=\operatorname{det}S^{\mathbf{J}_3}=\operatorname{det} 
\begin{blockarray}{ccccc}
  & 2 & 1 & 4 & 5 \\
\begin{block}{c(cccc)}
 A & -2 & 1 & 1 & 0 \\
 B & 2 & -1 & 0 & -1 \\
 C & 0 & 0 & 0 & -1 \\
 D & 0 & 0 & 0 & -1 \\
\end{block}
\end{blockarray}= 0.$$
Moreover, the two nonzero Child Selections $\mathbf{J}_1$ and $\mathbf{J}_2$ differ only in the image of species $A$: $\mathbf{J}_1(A)=1$, $\mathbf{J}_2(A)=2$. Consequently, by property of the determinant, $\alpha_{\mathbf{J}_1}= -2 \; \alpha_{\mathbf{J}_2}$. In particular:
$$\alpha_{\mathbf{J}_1}=\operatorname{det}
\begin{blockarray}{ccccc}
  & 1 & 3 & 4 & 5 \\
\begin{block}{c(cccc)}
 A & 1 & 0 & 1 & 0 \\
 B & -1 & 0 & 0 & -1 \\
 C & 0 & 1 & 0 & -1 \\
 D & 0 & 0 & 0 & -1 \\
\end{block}
\end{blockarray}
=1 \quad \text{and} \quad \alpha_{\mathbf{J}_2}=\operatorname{det}
\begin{blockarray}{ccccc}
  & 1 & 3 & 4 & 5 \\
\begin{block}{c(cccc)}
 A & -2 & 0 & 1 & 0 \\
 B & 2 & 0 & 0 & -1 \\
 C & 0 & 1 & 0 & -1 \\
 D & 0 & 0 & 0 & -1 \\
\end{block}
\end{blockarray}=-2.$$
Via Proposition \ref{vita}, the determinant of $G$ reads  
$$\operatorname{det}G=(\alpha_{\mathbf{J}_1}r_{1A}+\alpha_{\mathbf{J}_2}r_{2A})r_{3B}r_{4C}r_{5D}$$ 
and thus $\operatorname{det}G=0$ if, and only if, $r_{1A}=2r_{2A}$.\\ 

As explained in section \ref{almult}, the multiplicity of the eigenvalue zero can be asserted by looking at the trace of the adjugate matrix. We expands two diagonal entries, only, for sake of exemplification. Let us consider the diagonal entries $\operatorname{Adj}(G)^A_{A}$, and $\operatorname{Adj}(G)^D_D$. By Proposition \ref{AdG}, we have:
$$\operatorname{Adj}(G)^A_{A}=\sum_{\mathbf{J}^{\vee A}} \beta_{\mathbf{J}^{\vee A}}\prod_{m \neq A} r_{\mathbf{J}^{\vee A}(m)m}\quad\quad \text{and} \quad\quad \operatorname{Adj}(G)^D_{D}=\sum_{\mathbf{J}^{\vee D}}\beta_{\mathbf{J}^{\vee D}}\prod_{m \neq D} r_{\mathbf{J}^{\vee D}(m)m}.
$$
There are 2 Partial Child Selections $\mathbf{J}^{\vee A}$:
\begin{equation}
\begin{cases}
\mathbf{J}_1^{\vee A}:=\mathbf{J}_1^{\vee A}(B,C,D)=(3,4,5)\\
\mathbf{J}_2^{\vee A}:=\mathbf{J}_2^{\vee A}(B,C,D)=(1,4,5)\\
\end{cases}
\end{equation}
In particular, both $\mathbf{J}_1^{\vee A}$ and $\mathbf{J}_2^{\vee A}$ select reaction $4$. Note that 
$$S^{4}_{\vee A}=\begin{blockarray}{cc}
  & 4 \\
\begin{block}{c(c)}
 B & 0\\
 C & 0\\
 D & 0\\
\end{block}
\end{blockarray}\;,$$
hence both $\mathbf{J}_1^{\vee A}$ and $\mathbf{J}_2^{\vee A}$ are zero:
$$
\beta_{\mathbf{J}_1^{\vee A}}=\operatorname{det} 
\begin{blockarray}{cccc}
  & 3 & 4 & 5 \\
\begin{block}{c(ccc)}
 B & 0 & 0 & -1 \\
 C & 1 & 0 & -1\\
 D & 0 & 0 & -1\\
\end{block}
\end{blockarray}=0=\operatorname{det} 
\begin{blockarray}{cccc}
  & 3 & 4 & 5 \\
\begin{block}{c(ccc)}
 B & -1 & 0 & -1 \\
 C & 0 & 0 & -1\\
 D & 0 & 0 & -1\\
\end{block}
\end{blockarray}=\beta_{\mathbf{J}_2^{\vee A}}.$$
The example is designed so that a similar intuition implies that the diagonal entries $\operatorname{Adj}(G)^B_{B}$ and $\operatorname{Adj}(G)^C_{C}$ are zero, as well: we omit the analogous computation. With regard of $\operatorname{Adj}(G)^D_{D}$, there are 8 Partial Child Selections $\mathbf{J}^{\vee D}$:
$$
\begin{cases}
\mathbf{J}_1^{\vee D}:=\mathbf{J}_1^{\vee D}(A,B,C)=(2,1,4)\\
\mathbf{J}_2^{\vee D}:=\mathbf{J}_2^{\vee D}(A,B,C)=(2,1,5)\\
\mathbf{J}_3^{\vee D}:=\mathbf{J}_3^{\vee D}(A,B,C)=(1,3,4)\\
\mathbf{J}_4^{\vee D}:=\mathbf{J}_4^{\vee D}(A,B,C)=(2,3,4)\\
\mathbf{J}_5^{\vee D}:=\mathbf{J}_5^{\vee D}(A,B,C)=(1,5,4)\\
\mathbf{J}_6^{\vee D}:=\mathbf{J}_6^{\vee D}(A,B,C)=(2,5,4)\\
\mathbf{J}_7^{\vee D}:=\mathbf{J}_7^{\vee D}(A,B,C)=(1,3,5)\\
\mathbf{J}_8^{\vee D}:=\mathbf{J}_8^{\vee D}(A,B,C)=(2,3,5)\\
\end{cases}
$$
Analogously to the `full' Child Selection $\mathbf{J}_3$, both $\mathbf{J}_1^{\vee D}$ and $\mathbf{J}_2^{\vee D}$ are zero, as they select both reactions $1$ and $2$ whose stoichiometry is linearly dependent.
$$
\beta_{\mathbf{J}_1^{\vee D}}=\operatorname{det} 
\begin{blockarray}{cccc}
  & 2 & 1 & 4 \\
\begin{block}{c(ccc)}
 A & -2 & 1 & 1\\
 B & 2 & -1 & 0\\
 C & 0 & 0 & 0\\
\end{block}
\end{blockarray}=0=\operatorname{det} 
\begin{blockarray}{cccc}
  & 2 & 1 & 5 \\
\begin{block}{c(ccc)}
 A & -2 & 1 & 0\\
 B & 2 & -1 & -1\\
 C & 0 & 0 & -1\\
\end{block}
\end{blockarray}=\beta_{\mathbf{J}_2^{\vee D}}.$$
The other six Partial Child Selections can be grouped in three pairs: $(\mathbf{J}_3^{\vee D},\mathbf{J}_4^{\vee D})$, $(\mathbf{J}_5^{\vee D},\mathbf{J}_6^{\vee D})$, and $(\mathbf{J}_7^{\vee D},\mathbf{J}_8^{\vee D})$. Each pair consists of two Partial Child Selections which differ only in the image $\mathbf{J}^{\vee D}(A)$ of the species $A$: either the reaction $1$ or $2$.  Hence, the coefficients behaviors within each pair have a ratio of -2. In conclusion, again, in total analogy to the full Child Selections $\mathbf{J}_1$ and $\mathbf{J}_2$, we have that:
\begin{equation}
\begin{split}
\operatorname{Adj}(G)^D_{D}=&\sum_{\mathbf{J}^{\vee D}}\prod_{m \neq D} r_{\mathbf{J}^{\vee D}(m)m}\\
=&\sum_{i=3}^8 \beta_{\mathbf{J}_i^{\vee D}}\prod_{m \neq D} r_{\mathbf{J}_i^{\vee D}(m)m}\\
=&(\beta_{\mathbf{J}_3^{\vee D}}r_{1A}+\beta_{\mathbf{J}_4^{\vee D}}r_{2A})r_{3B}r_{4C} +(\beta_{\mathbf{J}_5^{\vee D}}r_{1A}+\beta_{\mathbf{J}_6^{\vee D}}r_{2A})r_{5B}r_{4C}\\
&+(\beta_{\mathbf{J}_7^{\vee D}}r_{1A}+\beta_{\mathbf{J}_8^{\vee D}}r_{2A})r_{3B}r_{5C}\\
=&(-r_{1A}+2r_{2A})r_{3B}r_{4C}+(r_{1A}-2r_{2A})r_{5B}r_{4C}+(r_{1A}-2r_{2A})r_{3B}r_{5C}\\
=&(r_{1A}-2r_{2A})(r_{5B}r_{4C}+r_{3B}r_{5C}-r_{3B}r_{4C}).
\end{split}
\end{equation}

Hence, $\operatorname{tr}\operatorname{Adj}(G)=(r_{1A}-2r_{2A})(r_{5B}r_{4C} +r_{3B}r_{5C} -r_{3B}r_{4C})$, and  $\operatorname{tr}\operatorname{Adj}(G)=0$ whenever $\lambda=(r_{1A}-2r_{2A})=0$, which characterizes $\operatorname{det}G=0$.

\section{Discussion}\label{discussion}

In this paper, we have presented an example of a chemical reaction network, for which the Jacobian of the associated dynamical system can be singular but never possesses a simple eigenvalue zero, for any choice of reaction rates $\mathbf{r}$ and any value of the concentrations $x$. The construction relies on studying algebraically the structure of the zero eigenvalues of the Jacobian $G$, using as a tool the adjugate matrix $\operatorname{Adj}(G)$.\\ 

It is natural to ask whether it is possible to have smaller examples of such a phenomenon. First, species $D$ and reactions $5$ and $F_D$ do not play a role in the algebraic feature we have presented. In particular, the same computation holds for a network $\tilde{\Gamma}$:
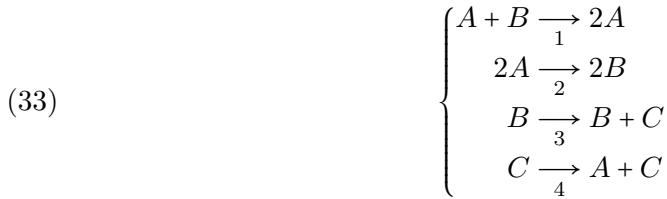
\begin{equation}
\begin{cases}
\begin{split}
A+B &\underset{1}{\longrightarrow} 2A\\
2A &\underset{2}{\longrightarrow} 2B\\
B &\underset{3}{\longrightarrow} B+C\\
C &\underset{4}{\longrightarrow} A+C\\
\end{split}
\end{cases}
\end{equation}
where species $D$ and reactions $5$ and $F_D$ have been removed. This constitutes the core example of the feature. However, the associated ODEs system, 
\begin{equation}
\begin{cases}
\dot{x}_A=r_1(x_A,x_B)-2 r_2(x_A) + r_4(x_C)\\
\dot{x}_B=-r_1(x_A,x_B)+2 r_2(x_A)\\
\dot{x}_C= r_3(x_B)\\
\end{cases},
\end{equation}
never admits an equilibrium, for any choice of $\mathbf{r}$. Given the intended application for equilibria bifurcation analysis, we have opted for an example that admits an equilibrium, at least. Nevertheless, let us stress one last time that the analysis is not linked with the precise value of the concentration $x$. It is opinion of the author that such feature cannot happen in networks with less than 3 species, even when $x$ is not necessarily an equilibrium.\\

The example we have presented includes reactions 1, 3, and 4, which are explicitly autocatalytic. Here, \emph{explicitly autocatalytic} simply means that species with nonzero stoichiometric coefficients appear at both sides of the reaction. It is possible to remove explicit autocatalysis by considering intermediates, hence splitting reactions 1, 3, 4 into
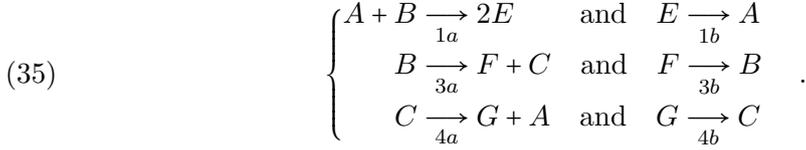
\begin{equation}
\begin{cases}
\begin{split}
A+B &\underset{1a}{\longrightarrow} 2E \quad\quad \; \text{and} \quad E \underset{1b}{\longrightarrow} A\\
B &\underset{3a}{\longrightarrow}  F + C \quad \text{and} \quad F \underset{3b}{\longrightarrow}B\\
C &\underset{4a}{\longrightarrow} G + A \quad \text{and} \quad G \underset{4b}{\longrightarrow} C\\
\end{split}
\end{cases}.
\end{equation} 
The computation follows analogously as the presented example. The removal of explicit autocatalysis by adding intermediates has the consequence of considerably increasing the size of the system. Thus, we have chosen to present the autocatalytic representation to keep the example as small as possible.\\

The mathematical literature on chemical reaction networks is often concerned with \emph{mass action} kinetics:
\begin{equation}
r_j(x)=k_j \prod_{m \in \mathbf{M}} x^{s^j_m}_m,
\end{equation}
where $s^j_m$ is the stoichiometric coefficient of $m$ as a reactant of the reaction $j$. Such an assumption gives rise to polynomial systems of differential equations with the interesting constraint that $x$ is real positive. Our system of ordinary differential equations \eqref{mainexsyst}, endowed with mass action, reads
\begin{equation}\label{mainma}
\begin{cases}
\dot{x}_A= k_1x_Ax_B-2k_2x_A^2 + k_4x_C\\
\dot{x}_B=-k_1x_Ax_B+2k_2x_A^2-k_5x_Bx_Cx_D\\
\dot{x}_C= k_3 x_B-k_5x_Bx_Cx_D\\
\dot{x}_D=-k_5x_Bx_Cx_D+F_D\\
\end{cases}.
\end{equation}
Unfortunately, this polynomial system does not admit an equilibrium with singular Jacobian, for any choice of $\mathbf{k}=(k_1,k_2,k_3,k_4,k_5,F_D)^T$ and $x>0$. In fact, the condition for a singular Jacobian is $r_{1A}=2r_{2A}$, as computed in section \ref{mainex}. According to mass action, the condition reads 
\begin{equation}\label{detma}
k_1 x_B =  4 k_2 x_A.
\end{equation}
Solving \eqref{mainma} for equilibria gives the additional constraint:
\begin{equation}\label{eqma}
2k_2x_A^2 - k_1x_Ax_B = k_5x_Bx_Cx_D > 0.
\end{equation}
But, inserting \eqref{detma} into \eqref{eqma} we obtain
\begin{equation}
2k_2x_A^2 - k_1x_Ax_B=2k_2x_A^2 - 4 k_2 x_A^2 < 0,
\end{equation}
hence there are no positive equilibra with singular Jacobian. However, it is easily possible to modify the system in order to maintain its validity as a mass action example. One possibility, with 6 species and 9 reactions is:
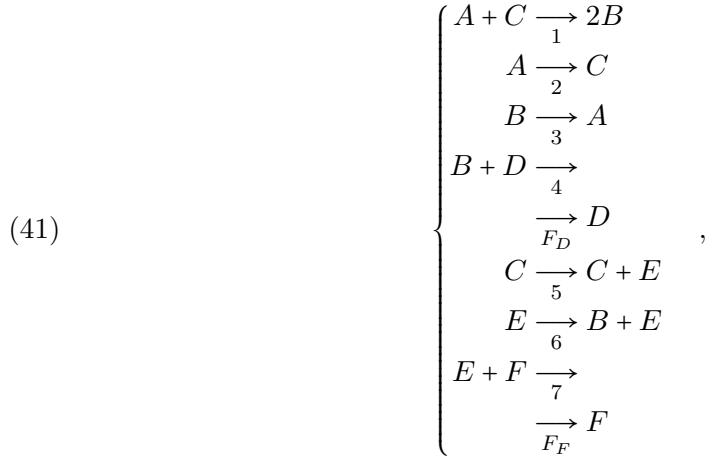
\begin{equation}
\begin{cases}
\begin{split}
A+C &\underset{1}{\longrightarrow} 2B\\
A &\underset{2}{\longrightarrow} C\\
B &\underset{3}{\longrightarrow} A\\
B+D &\underset{4}{\longrightarrow}\\
&\underset{F_D}{\longrightarrow}D\\
C &\underset{5}{\longrightarrow} C+E\\
E &\underset{6}{\longrightarrow} B+E\\
E+F &\underset{7}{\longrightarrow} \\
&\underset{F_F}{\longrightarrow} F\\
\end{split}
\end{cases},
\end{equation}
with associated mass action system
\begin{equation}\label{mainexsystma}
\begin{cases}
\dot{x}_A=-k_1x_Ax_C-k_2x_A +k_3 x_B \\
\dot{x}_B=2k_1x_A x_C -k_3 x_B-k_4x_Bx_D +k_6 x_E\\
\dot{x}_C=-k_1x_Ax_C+k_2x_A\\
\dot{x}_D=F_D -k_4x_Bx_D\\
\dot{x}_E=k_5x_C - k_7 x_E x_F\\
\dot{x}_F=F_F - k_7 x_E x_F\\
\end{cases}
\end{equation}
Again, we had to pay a price in terms of the dimension of the system. A straighforward computation shows: 
\begin{equation}
\begin{cases}
\operatorname{det}G=(k_1x_C-k_2)k_3k_4k_5 k_6 k_7 x_B x_E\\
\operatorname{tr}\operatorname{Adj}(G)=-(k_1x_C-k_2) k_3k_5 k_6 (k_4x_B + k_7 x_E)\\
\end{cases}.
\end{equation}
In particular, $\operatorname{det}G=0$ implies $\operatorname{tr}\operatorname{Adj}(G)=0$, hence system \eqref{amsystex} is never satisfied: it is possible to have a singular Jacobian but never with a simple eigenvalue zero, precisely as \eqref{mainexsyst}. Moreover, for $$\mathbf{k}=(k_1,k_2,k_3,k_4,k_5,k_6,k_7,F_D)^T=(1,1,2,1,1,1,1,1)^T,$$ the value 
$$x=(x_A, x_B, x_C, x_D, x_E, x_F)^T=(1,1,1,1,1,1)^T$$ is an equilibrium with a singular Jacobian. Since the core mathematical intuition is analogous, we have extendedly presented in terms of Child Selections and Partial Child Selections the more concise version \eqref{mainexsyst}, only.\\

The linearization at any zero-eigenvalue point of the presented system \eqref{mainexsyst} possesses a geometrically simple but algebraically double eigenvalue zero. This spectral condition is precisely the one of a \emph{Takens-Bogdanov bifurcation}. Such type of bifurcation was studied independently by Floris Takens \cite{Takens74} and Rifkat Bogdanov \cite{bog76}. In a neighborhood of the bifurcation point, it is possible to identify a saddle-node bifurcation curve, a Hopf bifurcation curve, and a homoclinic saddle connection curve. Hence, a Takens-Bogdanov bifurcation implies both multistationarity and oscillations due to the presence of, both, saddle-node and Hopf bifurcations, respectively. In a biochemical context, Kreusser and Rendall \cite{KrRen21} have proved the existence of a periodic orbit in a system modeling the activation of the \emph{Lymphocyte-specific protein tyrosine kinase}, by identifying a Takens-Bogdanov bifurcation. However, a Takens-Bogdanov bifurcation happens generically in a system with two (!) parameters, at an equilibrium whose Jacobian satisfies such spectral condition. On the contrary, in our example, the spectral condition is achieved by solving only one equality: $\operatorname{det}G=0$. A proper unfolding as Takens-Bogdanov is thus not possible. In this sense, our example can be considered a case of \emph{``not unfoldable Takens-Bogdanov bifurcation''}. The precise local dynamics of our example cannot be explained by standard bifurcation theorems, and it requires further investigation.\\

In conclusion, properties that are generic in general vector fields need not be generic when restricted to systems with a fixed network structure. A detailed investigation for further properties is needed. For instance, another theoretical possible scheme for bistability is a \emph{pitchfork bifurcation}. Pitchfork bifurcations relate to saddle-node in having the same spectral condition, but the bifurcation happens from a reference equilibrium persisting at any value of the bifurcation parameter, in contrast to the general saddle-node bifurcation. The equilibria diagram of a pitchfork bifurcation is topologically just the superposition of the diagram of a saddle-node bifurcation to a constant equilibrium line, suggesting the shape of a pitchfork. More in detail, at the bifurcation point, a stable equilibrium loses stability by generating two other stable equilibria (bistability!). On one side of the bifurcation point, there is one stable equilibrium, and on the other side, there are three equilibria: two stable equilibria with one unstable equilibrium within. A picture with reverse stability is of course possible, in analogy. In contrast to saddle-node bifurcations, however, pitchfork bifurcations are nongeneric in the set of vector fields with a singular Jacobian, hence - as already discussed - in applications we expect bistability arising from a pair of connected saddle-node bifurcations, rather than from one single pitchfork bifurcation. Investigating whether certain networks have the property to exhibit generically pitchfork bifurcations rather than saddle-node bifurcations is particularly interesting as it provides an alternative and unexpected scheme for bistability, where no hysteretic switch-like behavior is present. More in general, this work calls for investigating further whether the network structure may or may not interfere with the genericity of certain properties. Such a question can be addressed both in the case of a chosen kinetics, or more in general as in this paper. Answers are quite interesting both ways: a positive answer would greatly simplify the bifurcation analysis, as it would grant for free such generic property, while a negative answer would provide examples of networks that show surprising features.\\

\bibliography{bibliography/references.bib}
 \bibliographystyle{ieeetr}
 
\end{document}